\documentclass[10pt]{article}
\usepackage[tbtags]{amsmath}
\usepackage{epsfig,amstext,amssymb,amsthm,latexsym}
\usepackage{amssymb,color}
\definecolor{c20}{rgb}{0.,0.7,0.}
\definecolor{c30}{rgb}{0.,0.,1.}
\definecolor{c40}{rgb}{1,0.1,0.7}
\definecolor{c50}{rgb}{1,0,0}

\def\pE#1{\textcolor{c30}{#1}}
\def\pE#1{#1}
\def\cE#1{\textcolor{c30}{#1}}
\def\cE#1{#1}
\def\cQ#1{\textcolor{c30}{#1}}
\def\pE#1{#1}
\def\cQ#1{#1}

\def\ZW#1{\textcolor{c20}{#1}}
\def\ZW#1{#1}

\def\cW#1{{\textcolor{c40}{#1}}}
\def\cW#1{#1}

\def\CZ#1{{\textcolor{c40}{#1}}}
\def\CZ#1{#1}
\def\CCZ#1{{\textcolor{c40}{#1}}}
\def\CCZ#1{#1}
\def\CCC#1{{\textcolor{c40}{#1}}}
\def\CCC#1{#1}
\def\CZC#1{{\textcolor{c40}{#1}}}
\def\CZC#1{#1}
\def\czw#1{{\textcolor{c40}{#1}}}
\def\czw#1{#1}

\catcode`\@=11 \@addtoreset{equation}{section} \catcode`\@=12

\allowdisplaybreaks[4] 
\newcommand{\nwc}{\newcommand}
\nwc{\COM}[1]{}

\nwc{\vs}[1]{\vskip #1 cm}

\newtheorem{theo}{Theorem}[section]
\newtheorem{sat}[theo]{Proposition}
\newtheorem{de}[theo]{Definition}
\newtheorem{lem}[theo]{Lemma}

\newtheorem{korr}[theo]{Corollary}
\newtheorem{remark}[theo]{Remark}
\newtheorem{exxa}[theo]{Example}

\newcommand{\netheo}[1]{{Theorem \ref{#1}}}
\newcommand{\nekorr}[1]{{Corollary \ref{#1}}}

\def\d{\mathrm{d}}

\newcommand{\ve}{\varepsilon}

\newcommand{\abs}[1]{\lvert #1 \rvert}
\newcommand{\Abs}[1]{ \Bigl \lvert #1 \Bigr \rvert}

\newcommand{\pb}[1]{\mathbb{P}\Bigl \{#1 \Bigr \}}
\newcommand{\pk}[1]{\mathbb{P}\left\{#1 \right\}}
\newcommand{\R}{\!I\!\!R}

\newcommand{\inr}{\in \R}

\newcommand{\limit}[1]{\lim_{#1 \to   \infty}}

\newcommand{\equaldis}{\stackrel{d}{=}}

\newcommand{\BQN}{\begin{eqnarray}}
\newcommand{\EQN}{\end{eqnarray}}
\newcommand{\BQNY}{\begin{eqnarray*}}
\newcommand{\EQNY}{\end{eqnarray*}}
\newcommand{\BS}{\begin{sat}}
\newcommand{\ES}{\end{sat}}
\newcommand{\BL}{\begin{lem}}
\newcommand{\EL}{\end{lem}}
\newcommand{\BT}{\begin{theo}}
\newcommand{\ET}{\end{theo}}
\newcommand{\BK}{\begin{korr}}
\newcommand{\EK}{\end{korr}}
\newcommand{\BD}{\begin{de}}
\newcommand{\ED}{\end{de}}
\newcommand{\BIT}{\begin{itemize}}
\newcommand{\EIT}{\end{itemize}}
\newcommand{\BDI}{\begin{description}}
\newcommand{\EDI}{\end{description}}

\newcommand{\BEX}{\begin{exxa}}
\newcommand{\EEX}{\end{exxa}}

\newcommand{\QED}{\hfill $\Box$}
\newcommand{\IF}{\infty}

\def\fracl#1#2{\biggr( \frac{#1}{#2} \biggl) }
\newcommand{\EXP}[1]{\exp \left( #1 \right) }

\newcommand{\prooftheo}[1]{ \textsc{Proof of Theorem} \ref{#1} }

\newcommand{\proofkorr}[1]{\textsc{Proof of Corollary} \ref{#1}}

\begin{document}

\centerline{\bf \large Tail Asymptotic of Weibull-Type Risks}

  \vskip 1.5 cm
\centerline{Enkelejd Hashorva and Zhichao Weng\thanks{
Research supported by the Swiss National Science Foundation Project 200021-134785} }

  \vskip .1 cm
\centerline{Department of Actuarial Science, University of Lausanne, Switzerland}

  \vskip .5 cm
\centerline{  \today{}}

  \vskip .1 cm
{\bf Abstract:} With motivation from Arendarczyk  and D\c{e}bicki (2011)\CZC{,} in this paper we derive the tail asymptotics of the product of two dependent \CCC{Weibull}-type risks\CZC{,} which is of interest in various statistical and applied probability problems.  Our \pE{results extend} some
\pE{recent findings of Schlueter and Fischer (2012)} and Bose et al.\ (2012).


\section{Introduction}

In numerous statistical and probabilistic models various quantities of interest are defined in terms of product of random variables (or risks).
For instance, given $X_1,X_2$ two positive risks,  the product $Z= X_1X_2$ can be used to model a random deflation/inflation effect,
if say $X_1$ is the deflator/inflator and $X_2$ is some base risk related to some financial loss.
Since often the distribution \cE{functions} of the risks are not known, the main interest is on the asymptotic analysis of the tail of $Z$.
When $X_1$ is a bounded random variable, then $Z$ can be seen as a random contraction of $X_2$, see e.g., Berman (1992), Cline and Samorodnitsky (1994), Pakes and Navarro (2007), Hashorva and Pakes (2010),  Hashorva et al.\ (2010,\czw{2012}), Hashorva (\czw{2011},2012), Yang and Wang (2012). Interesting models where $X_1$ is unbounded have been studied in Cline and Samorodnitsky (1994), Maulik and Resnick (2004), Nadarajah (2005), Nadarajah and Kotz (2005), Zwart et al.\ (2005), Jessen and Mikosch (2006), Tang (2006a,b,2008), Liu and Tang (2010), Arendarczyk  and D\c{e}bicki (2011,2012), Balakrishnan and Hashorva (2011), Chen (2011), Constantinescu et al.\ (2011),  Jiang and Tang (2011), Yang et al.\ (2011), Schlueter and Fischer (2012) among several others.

With motivation from Arendarczyk and D\c{e}bicki (2011),  in this paper we are concerned with the investigation of the tail asymptotics of the product $Z=X_1X_2$ of risks with \CZC{Weibull} tail behaviour i.e., for $X_1$ and $X_2$ such that
\BQN\label{Weib}
  \pk{X_i> x} \sim g_i(x) \exp(- L_i x^{p_i}),
\EQN
with $g_i(\cdot)$ some regularly varying function at infinity and $L_i,p_i,i=1,2$ positive constants.
In our notation $a(x)\sim b(x)$\CCC{,} for two functions \cE{$a(\cdot)$ and $b(\cdot)$}\CCC{,} means that $\limit{x}a(x)/b(x)=1$.

A large class of such risks satisfy \eqref{Weib} with  $g_i(\cdot)$ a polynomial function i.e.,
\begin{eqnarray}\label{poly}
  \pk{X_i> x} \sim C_i x^{\alpha_i} \exp(- L_i x^{p_i}),
\end{eqnarray}
 with $C_i,p_i,L_i,i=1,2$ positive constants, $\alpha_1,\alpha_2\in \R$. A remarkable result of Arendarczyk and D\c{e}bicki (2011)\CCC{,} which is crucial for the analysis of the extremes of Gaussian processes over random intervals\CCC{,} shows that when \eqref{poly} holds, then
\begin{eqnarray}\label{DE}
\pk{Z> x}\sim\Bigl(\frac{2\pi p_{2}L_{2}}{p_{1}+p_{2}}\Bigr)%
^{\CZC{\frac{1}{2}}}C_{1}C_{2}A^{\CZC{\frac{p_{2}}{2}}+\alpha _{2}-\alpha _{1}}x^{\frac{2p_{2}\alpha
_{1}+2p_{1}\alpha _{2}+p_{1}p_{2}}{2(p_{1}+p_{2})}}
 \exp \left( -Bx^{\frac{p_{1}p_{2}}{%
p_{1}+p_{2}}}\right),
\end{eqnarray}
with
 \BQN\label{CA}
 A =[(p_1L_1)/(p_2 L_2)]^{1/(p_1+p_2)} \mbox{   and   } B=L_{1}A^{-p_{1}}+L_{2}A^{p_{2}}.
 \EQN
Clearly, also $Z$ is Weibull-type risk, and thus \eqref{DE} shows the closure property \cE{for the product of such risks.}

The main goal of this paper is to investigate the tail asymptotics of $Z$ for \CCC{Weibull}-type risks 
allowing further
for the risks to be dependent. In various theoretical problems and applications independence assumption is not tenable.
Particular examples of the dependence structure assumed in this paper are risks with bivariate
Fairly-Gumbel-Morgenstern \CZC{(FGM)} distribution. We \CZC{also} show by considering the special case that $X_1$ and $X_2$ are jointly Gaussian,
that the dependence structure is crucial for the tail asymptotic of $Z$. \\
 Our findings are of interest in various probabilistic models, for instance our \nekorr{co1} subsumes Theorem 1 in Bose et al. (2012)
which \CCC{is} crucial for dealing with the spectral radius of random k-circulants; in particular that result implies \cE{the} closure property of independent Weibull-type risks with respect to product.
\cE{Our first application deals with the supremum of Brownian motion over random time interval.
 In the second application we extend the findings of Schlueter and Fischer (2012) which concern the calculation of the weak tail dependence coefficient of elliptical generalized hyperbolic distribution.}

Outline of the rest of the paper: Section 2 presents the main findings of this contribution. In Section 3 we give two applications, followed by Section 4 where \cE{all} the proofs are displayed. 

\section{Main Results}
\cE{In this section both risks $X_1 \sim F_1$ and $X_2 \sim F_2$ are positive and satisfy}
\eqref{Weib}  with $p_i,L_i$ positive constants, and $g_i$ regularly varying at infinity with index $\alpha_i, i=1,2$.
Their dependence structure is modeled by a tractable conditions, namely we shall assume that
for some positive measurable function $c(\cdot, \cdot)$ and some constants $K_1>0$, $K_2>0$, $\beta_1, \beta_2 \inr $
\BQN\label{depB}
\pk{X_1 > x/y \lvert X_2=y} =\pk{X_1> x/y} c(x, y) \quad \text{ and } \quad K_1 x^{\beta_1} \le c(x,y) \le  K_2 x^{\beta_2}
\EQN
\CCC{are} satisfied for all $x$ large and any $y>0$ and further
\BQN\label{depB2}
\lim_{x\to \IF} \sup_{  \ZW{y}\in [a_1 w_x, a_2 w_x]} \Abs{ c(x,y)- D x^{q_1} \ZW{y^{q_2-q_1}}} =0
\EQN
holds for some constants $D>0, 0 < a_1< a_2$, $q_1,q_2\inr$ and $w_x=x^{\frac{p_1}{p_1+p_2}}$.\\


\BT\label{th1}
\cE{Let $X_1$ and $X_2$ be two dependent risks as above such that both $g_1,g_2$ are ultimately monotone}. If condition \eqref{depB} and \eqref{depB2} hold, then
\BQNY
\pk{Z>x}\sim
D \Bigl(\frac{2\pi p_{2}L_{2}}{p_{1}+p_{2}}\Bigr)
^{\CZC{\frac{1}{2}}}A^{\CZC{\frac{p_2}{2}}+q _{2}-q _{1}}x^{\frac{2p_{2} q
_{1}+2p_{1}q _{2}+p_{1}p_{2}}{2(p_{1}+p_{2})}}g_1(z_x^{-1}x)g_2(z_x)
\exp\left(-Bx^{\frac{p_1p_2}{p_1+p_2}}\right), 
\EQNY
with \CZ{$z_x=Ax^{p_1/(p_1+p_2)}$} and $A$ and $B$ given by \eqref{CA}.
\ET

\BK
\label{co1}
Under the conditions of Theorem \ref{th1}, and $X_1, X_2$ are independence, \CZC{then}
\begin{eqnarray}
\pk{Z>x} &\sim&\left(\frac{2\pi p_2L_2}{p_1+p_2}\right)^{\CZC{\frac{1}{2}}}A^{\frac{p_2}{2}}x^{\frac{p_1p_2}{2(p_1+p_2)}}
g_1(z_x^{-1}x)g_2(z_x)
\EXP{-Bx^{\frac{p_1p_2}{p_1+p_2}}}.
\end{eqnarray}
If additionally $X_1$ possess a positive pdf $h_1$ which is \CCZ{bounded and} ultimately decreasing, then the pdf $h$ of $Z$ satisfies
\BQN\label{pdf}
h(x)&\sim&L_1p_1A^{-p_1}x^{\frac{p_1p_2}{p_1+p_2}-1}\pk{Z>x},
\EQN
provided that $h_1(x)=(1+o(1))L_1p_1x^{p_1-1}\pk{X_1>x}$.
\EK

\cE{Bose et al. (2012) derived in their Theorem 1 the tail asymptotics of the product of \CZ{$n$} independent unit exponential random variables.
The above corollary extends Theorem 1 of the aforementioned paper to the product of independent Weibull-type risks with common parameters $L$ \CCZ{and $p$} and $g$ being ultimately monotone.  In fact, if  $X_i \sim F, i=1,\cdots, m$ \cE{are} independent positive random variables such that \eqref{poly} holds with $C, p, L$ positive constants, $\alpha \inr$, then Theorem 1 of the aforementioned paper can be generalised to the following statement
\BQN\label{DE2}
\pk{\prod^m_{i=1}X_i>x} \sim m^{-\frac{1}{2}}(2 \pi L)^{\frac{m-1}{2}}C^m x^{\frac{2m\alpha+(m-1)p}{2m}}\EXP{-mLx^{\frac{p}{m}}}, 
\EQN
which is a direct implication of the result of \eqref{DE} derived in Arendarczyk and D\c{e}bicki (2011).}

{\bf Remarks:} a) Liu and Tang (2010) considers more general  Weibull-type risks and establishes under weaker conditions than ours
the subexponentiality of $Z$. \\
b) If $K_1=0$ in condition \eqref{depB}, the lower bound of $\pk{Z>x}$ can be substituted by
\BQN\label{depB1}
\pk{X_1 > x, X_2> y} \ge K x^{\gamma_1}y^{\gamma_2}\pk{X_1> x} \pk{X_2> y}
\EQN
for all $x,y$ large and some constants $K>0$, $\gamma_1, \gamma_2 \inr$. \\
%
c) \czw{As can} be seen from the proof of \netheo{th1} (check in particular \eqref{uf}), the assumption that $g_i,i=1,2$ is regularly varying
can be slightly weakened to
\BQN
\lim_{\ve \to 0} \lim_{u \to \IF} \frac{g_i((1+ \ve) u)}{g_i(u)}=1, \quad
\text{ and   }
c_iu^{r_i}\le g_i(u) \le c_i^*u^{r_i^*},
\EQN
where the inequalities holds for all large $u$ with constants $c_i>0$, $r_i, r_i^* \inr$, $i=1,2.$\\
d) The constants appearing in the conditions (2.1) and (2.2) do not explicitly show in the tail asymptotics of $Z$.
Our dependence model implied by the aforementioned conditions is quite restrictive. As shown below in Example 3, complete different results are obtained if we drop some restrictions on the joint dependence of $X_1$ and $X_2$.

\COM{
\begin{remark} Particularly, as $x \to \infty$, we get first order expansion as following,
\begin{eqnarray*}
\pk{X_1X_2> x}\sim\Bigl(\frac{2\pi p_{2}L_{2}}{p_{1}+p_{2}}\Bigr)%
^{1/2}C_{1}C_{2}A^{p_{2}/2+\alpha _{2}-\alpha _{1}}x^{\frac{2p_{2}\alpha
_{1}+2p_{1}\alpha _{2}+p_{1}p_{2}}{2(p_{1}+p_{2})}}
 \exp ( -(L_{1}A^{-p_{1}}+L_{2}A^{p_{2}})x^{\frac{p_{1}p_{2}}{%
p_{1}+p_{2}}}),
\end{eqnarray*}
and second order expansion as following,
\begin{eqnarray*}
\pk{X_1X_2> x}
&\sim&\sqrt{2\pi }
C_{1}C_{2}p_{2}L_{2}A^{p_{2}+\alpha _{2}-\alpha _{1}}
x^{\frac{2p_{2}\alpha _{1}+2p_{1}\alpha _{2}+p_{1}p_{2}}{2(p_{1}+p_{2})}}\exp
(-\psi (1)x^{\frac{p_{1}p_{2}}{p_{1}+p_{2}}})\\
&&\times \Bigl(\frac{1}{\sqrt{\psi ^{\prime
\prime }(1)}}+
\left(\frac{5}{24}\frac{(\psi^{\prime \prime\prime}(1))^2}{(\psi^{\prime\prime}(1))^{\frac{7}{2}}}-
\frac{1}{8}\frac{\psi^{(4)}(1)}{(\psi^{\prime\prime}(1))^{\frac{5}{2}}}\right)
x^{-\frac{p_{1}p_{2}}{(p_{1}+p_{2})}}\Bigr)\\
&=&\Bigl(\frac{2\pi p_{2}L_{2}}{p_{1}+p_{2}}\Bigr)
^{1/2}C_{1}C_{2}A^{p_{2}/2+\alpha _{2}-\alpha _{1}}x^{\frac{2p_{2}\alpha
_{1}+2p_{1}\alpha _{2}+p_{1}p_{2}}{2(p_{1}+p_{2})}}
\exp ( -(L_{1}A^{-p_{1}}+L_{2}A^{p_{2}})x^{\frac{p_{1}p_{2}}{
p_{1}+p_{2}}}) \\
&&\times \Bigl(1+\frac{2(p_2-p_1-3)^2-3p_1p_2-6}{24L_2A^{p_2}p_2(p_1+p_2)}x^{-\frac{p_1p_2}{p_1+p_2}}
\Bigr).
\end{eqnarray*}
c) If in \nekorr{co2} $\alpha=0,L=1=p=1$ then we obtain an explicit form of the claim of Theorem 1 of Bose et al. (2012).
\end{remark}
}
\COM{

\BT
\label{th3}
Let $X_{i}, i=1,2$ be two positive independent random variables such that
\begin{eqnarray*}
  \pk{X_i> x} \sim C_i x^{\alpha_i} \exp(- L_i x^{p_i}),
\end{eqnarray*}
 with $C_i,p_i,L_i,i=1,2$ positive constants, $\alpha_1,\alpha_2\in \R$.
 If further $X_2$ possess a positive pdf  $h_2$ such that $h_2(x) =(1+o(1))\pk{X_2> x}  L_2p_2 x^{p_2-1}$,
then with $ A =   [(p_1L_1)/(p_2 L_2)]^{1/(p_1+p_2)}$ as $x\rightarrow \infty $ we have
\begin{eqnarray*}
\pk{X_1X_2> x} &=&(1+o(1))
\Bigl(\frac{2\pi p_{2}L_{2}}{p_{1}+p_{2}}\Bigr)
^{1/2}C_{1}C_{2}A^{p_{2}/2+\alpha _{2}-\alpha _{1}}x^{\frac{2p_{2}\alpha
_{1}+2p_{1}\alpha _{2}+p_{1}p_{2}}{2(p_{1}+p_{2})}}\\
&&\exp(-(L_1A^{-p_1}+L_2A^{p_2})x^{\frac{p_1p_2}{p_1+p_2}})
\left(1+\sum^{\infty}_{n=1}\lambda_n x^{-\frac{np_{1}p_{2}}{(p_{1}+p_{2})}}\right)
\end{eqnarray*}
where the coefficients $\lambda_n$ are polynomials of the higher derivatives $\psi^{(k)}(1), k\geq2$ with
 $\psi
(y)=L_{1}A^{-p_{1}}y^{-p_{1}}+L_{2}A^{p_{2}}y^{p_{2}}$.
\ET
}

\COM{
\textbf{Example 1.} Let $X^2_i,i=1,2,\cdot\cdot\cdot,m$ be chi-squared distributed with
$k$ degrees of freedom i.e.,
\begin{eqnarray*}
\pk{X_i>x}\sim \frac{x^{k-2}}{2^{k/2-1}\Gamma(k/2)}\exp(-x^2/2)
\end{eqnarray*}
as $x \to \infty$. Using Corollary \ref{co2}, we have
\begin{eqnarray*}
\Pk{\prod_{i=1}^{m}X_i> x} &\sim&
\left(\frac{\pi^{m-1}}{m}\right)^{\frac{1}{2}}\frac{1}{2^{mk/2-m}(\Gamma(k/2))^m}
x^{\frac{mk-m-1}{m}}\exp\left(-\frac{m}{2}x^{\frac{2}{m}}\right)\left(1+\sum^{\infty}_{n=1}\lambda_nx^{-\frac{2n}{m}}\right)
\end{eqnarray*}
where the coefficients $\lambda_n$ are polynomials of the higher derivatives $\psi^{(k)}(1), k\geq2$ with
 $\psi
(y)=\frac{1}{2}y^{-2}+\frac{m-1}{2}y^{\frac{2}{m-1}}$.
}

\cE{We present next three examples.}
\COM{
\textbf{Example 1.}
Let $X^2_i,i=1,2,\cdot\cdot\cdot,m$ be Gamma distributed with
scale $\lambda$ and shape $\alpha$ i.e.,
\begin{eqnarray*}
\pk{X_i>x}\sim \frac{x^{2\alpha-2}}{\lambda^{\alpha-1}\Gamma(\alpha)}\exp\left(-\CZC{x^2}{\lambda}\right).
\end{eqnarray*}
\cE{In view of} Corollary \ref{DE2}, we have
\begin{eqnarray*}
\pk{\prod_{i=1}^{m}X_i> x} &\sim&
\left(\frac{2^{m-1}\pi^{m-1}}{m\lambda^{m-1}}\right)^{\frac{1}{2}}\frac{1}{\lambda^{m\alpha-m}(\Gamma(\alpha))^m}
x^{\frac{2m\alpha-m-1}{m}}\exp\left(-\frac{m}{\lambda}x^{\frac{2}{m}}\right).
\end{eqnarray*}
\cE{Furthermore}, {by Corollary \ref{co1}, for the pdf $h$ of $\prod_{i=1}^{m}X_i$ we obtain}
$$h(x) \sim \frac{2}{\lambda}x^{\frac{2}{m}-1}\pk{\prod_{i=1}^{m}X_i> x}.$$}

\textbf{Example 1.}
Let $\CCZ{X_i},i=1,2,\cdot\cdot\cdot,m$ be Gamma distributed with
scale $\lambda$ and shape $\alpha$ i.e.,
\begin{eqnarray*}
\pk{X_i>x}\sim \frac{\CCZ{x^{\alpha-1}}}{\lambda^{\alpha-1}\Gamma(\alpha)}\exp(-\CCZ{x}/\lambda)
\end{eqnarray*}
as $x \to \infty$. \cE{In view of} \CCZ{\eqref{DE2}}, we have
\begin{eqnarray*}
\pk{\prod_{i=1}^{m}X_i> x} &\sim&
\left(\frac{2^{m-1}\pi^{m-1}}{m\lambda^{m-1}}\right)^{\frac{1}{2}}\frac{1}{\lambda^{m\alpha-m}(\Gamma(\alpha))^m}
x^{\frac{2m\alpha-m-1}{\CCZ{2m}}}\exp\left(-\frac{m}{\lambda}x^{\frac{\CCZ{1}}{m}}\right).
\end{eqnarray*}
\cE{Furthermore}, {by \CCZ{\eqref{pdf}}, we get for the pdf $h$ of $\prod_{i=1}^{m}X_i$}
$$h(x) \sim \frac{x^{\frac{\CCZ{1}}{m}-1}}{\lambda}\pk{\prod_{i=1}^{m}X_i> x}.$$

\textbf{Example 2.}
Let $X_i \sim F_i, i=1,2$ be two positive random variables such that \eqref{Weib} holds \cE{and $g_1,g_2$ are} ultimately monotone and regularly varying at infinity. \cE{We suppose that the joint distribution of $X_1$ and $X_2$ is FGM} i.e., for $\tau \in [-1,1]$,
\BQNY
\pk{X_1\leq z_1,X_2\leq z_2}=F_1(z_1)F_2(z_2)(1-\tau(1-F_1(z_1))(1-F_2(z_2))).
\EQNY
Consequently,
$$\pk{X_1>x/y|X_2=y}=\overline{F}_1(x/y)(1+\tau F_1(x/y)(1-2F_2(y)))$$
where
$$\CCZ{1-|\tau|}\leq 1+\tau F_1(x/y)(1-2F_2(y))<\CCZ{1+|\tau|},$$
\CCZ{and}
$$\pk{X_1>z_1,X_2>z_2}=\overline{F}_1(z_1)\overline{F}_2(z_2)(1-\tau F_1(z_1)F_2(z_2))\geq (1-\CCZ{|\tau|})\overline{F}_1(z_1)\overline{F}_2(z_2).$$
\cE{Hence both assumptions \eqref{depB} and \eqref{depB1} are satisfied for FGM dependence. Further,}
\BQNY
\lim_{x\to \IF} \sup_{  y\in [a w_x, a^{-1} w_x]}
\Abs{(1+\tau F_1(x/y)(1-2F_2(y))) - (1-\tau)} =0,
\EQNY
\cE{hence \CCC{the} condition \eqref{depB2} holds with $D=1- \tau$.}
\cE{A direct application of  Theorem \ref{th1}  yields}
\BQN \label{eq FGM}
\pk{Z>x}\sim
(1-\tau) A^{\frac{p_{2}}{2}}
\Bigl(\frac{2\pi p_{2}L_{2}}{p_{1}+p_{2}}\Bigr)^{\frac{1}{2}}
x^{\frac{p_{1}p_{2}}{2(p_{1}+p_{2})}}g_1(z_x^{-1}x)g_2(z_x)
\exp\left(-Bx^{\frac{p_1p_2}{p_1+p_2}}\right).
\EQN

\textbf{Example 3.} Let $X_i$, $i=1,2$ be two standard Gaussian random variables with correlation coefficient $\rho\cE{\in (-1,1)}$.
For this example the dependence function is different from that of FGM treated above. In particular condition \eqref{depB} is not satisfied
\CZ{since the conditional distributions are Gaussian}. After some straightforward calculations we  obtain 
\BQN\label{Gauss}
\pk{Z>x}\sim \pE{\frac{1+\rho}{\sqrt{2\pi x}}}\EXP{-\frac{x}{1+\rho}}.
\EQN
\pE{Note that when $\rho=0$, then \eqref{Gauss} follows directly by \CCC{\eqref{DE}}.}
\cE{The asymptotics in \eqref{Gauss} shows that instead of $B$ appearing in \eqref{eq FGM}, the term $1/(1+ \rho)$ which depends on the correlation coefficient $\rho$ appears. Our dependence structure does not imply restrictions for $B$, hence the Gaussian case is clearly not covered by the dependence model assumed in Theorem \ref{th1}.}

\section{Applications}
Our first application \cE{deals} with the supremum of Brownian motion on some random interval $[0,\mathcal{T}]$. It can be easily seen that our result can be extended for several Gaussian processes using the key findings of Arendarczyk  and D\c{e}bicki (2011).

Assume that $\mathcal{T}$ is almost surely positive with asymptotic tail behaviour given by \eqref{Weib} with some function $g(\cdot)$ and
positive \CCC{constants} $L,p$. If $B(t),t\ge 0$ is a standard Brownian motion (mean 0, variance function $t$, and continuous sample paths), then
for any $x>0$, by the self-similarity property of Brownian motion we have
\BQNY
\pb{\sup_{t \in [0,\mathcal{T}]} B(t)> x} = \pb{\mathcal{T}^{1/2} \sup_{t\in [0,1]} B(t)> x}\ge \pb{B(\mathcal{T})>x}.
\EQNY
Since $\sup_{t\in [0,1]} B(t)$ has the same distribution as $\abs{B(1)}$, if further $g(\cdot)$ is ultimately monotone,
applying Corollary \ref{co1} we obtain
\CZ{
\BQN
\pb{\sup_{t \in [0,\mathcal{T}]} B(t)> x} \sim \left(\frac{2}{1+p}\right)^{\frac{1}{2}}g\left(Ax^{\frac{2}{1+p}}\right)\EXP{-\left(\frac{1}{2A}+LA^{p}\right)x^{\frac{2p}{1+p}}},
\EQN
where \cQ{$A=(2Lp)^{-1/(1+p)}$}.
For the special simple case }$g(x)= C x^\alpha$ with $C$ some positive constant \cE{the above claim is} stated in Theorem 4.1 
of  Arendarczyk  and D\c{e}bicki (2011).

Our second application is motivated by the recent paper  Schlueter and Fischer (2012) which derives a formula for
the weak tail dependence coefficient of elliptical generalized hyperbolic distribution (EGHD).\\
We shall consider below a bivariate elliptical random vector $(X_1,X_2)$ with stochastic representation
\BQN\label{Ellip}
(X_1,X_2)\equaldis R(U_1,\rho U_1+\sqrt{1-\rho^2}U_2), \qquad \qquad \rho \in(-1,1),
\EQN
where the positive random radius $R$ is independent of $(U_1,U_2)$ which is uniformly
distributed on the unit circle of $\R^2$. The basic properties of elliptical random vectors are
well-known, see e.g., Cambanis et al. (1981). 
Assume that the random radius $R$ has distribution function $G$ in the Gumbel max-domain of attraction \pE{(see e.g., Resnick (1987))} i.e.,
\BQNY
\lim_{x \to \infty}\frac{1-G(x+s/ w(x))}{1-G(x)}=\EXP{-s}, \qquad \forall s\inr
\EQNY
holds with some positive scaling function $w$, Hashorva (2007) obtained the
exact asymptotic of tail probability of the bivariate elliptical vector
\BQN\label{elliptical}
\pk{X_1>x,X_2>x}\sim \frac{\czw{\sqrt{c_{\rho}}}}{2 \pi}\frac{(1-\rho^2)^{\frac{3}{2}}}{(1-\rho)^2}\frac{1}{x w(\sqrt{c_{\rho}}x)}
\pk{R>\sqrt{c_{\rho}}x},
\EQN
where $c_{\rho}=2/(1+\rho)$.\\
For statistical modelling, calculation of the weak tail dependence coefficient is of particular interest.
Hashorva (2010) derived the weak tail dependence coefficient of the elliptical distribution as
$$\chi=2\left(\frac{1+\rho}{2}\right)^{\theta/2}-1,$$
if
\BQNY
\lim_{x \to \infty}\frac{ w(cx)}{ w(x)}=c^{\theta-1}, \qquad \forall c>0,
\EQNY
holds for some $\theta \in [0, \infty)$.
We extend the above results to bivariate scaled elliptical random vectors under the condition that
the joint distribution of the random radius $R$ and the scaling random variable $S$ is the FGM distribution.\\

\BT\label{ap1}
Let $(X_1,X_2)$ be a bivariate elliptical random vector with representation  \eqref{Ellip} \cE{and define $Y_1=SX_1,Y_2=SX_2$ with $S$ some positive scaling random variable. Assume that both $R$ and }
$S$ \cE{satisfy \eqref{Weib} with $g_1,g_2$ ultimately monotone}, and have FGM distribution.
\CZ{If $SR$ is independent of $(U_1,U_2)$}, then we have
\BQNY
\pk{Y_1>x,Y_2>x}&\sim&
(1-\tau)\frac{(1-\rho^2)^{\frac{3}{2}}}{(1-\rho)^2}
\Bigl(\frac{p_{2}L_{2}}{2\pi (p_{1}+p_{2})}\Bigr)
^{\frac{1}{2}}c_{\rho}^{1-\frac{p_1p_2}{4(p_1+p_2)}}
(p_1L_1)^{-1}A^{\frac{p_{2}}{2}+p_1}x^{-\frac{p_{1}p_{2}}{2(p_{1}+p_{2})}}\\
&& \times
g_1\left(c_{\rho}^{\frac{p_2}{2(p_1+p_2)}}z_x^{-1}x\right)
g_2\left(c_{\rho}^{\frac{p_1}{2(p_1+p_2)}}z_x\right)
\exp\left(-Bc_{\rho}^{\frac{p_1p_2}{2(p_1+p_2)}}x^{\frac{p_1p_2}{p_1+p_2}}\right),
\EQNY
and the weak tail dependence coefficient of the random pair $(Y_1,Y_2)$ is given by
$$\chi=2\cdot\left(\frac{1+\rho}{2}\right)^{\frac{p_1p_2}{2(p_1+p_2)}}-1.$$
\ET

\textbf{Example 4.}
A canonical example of a bivariate scaled elliptical distribution is the EGHD,
which is now widely used in finance (see e.g., Eberlein and Keller (1995) and McNeil et al. (2005)).\\
\cE{Let $\CCZ{(Y_1,Y_2)}$} be elliptical generalized hyperbolic random vector with stochastic representation
$$(Y_1,Y_2)\equaldis (\CCZ{S}X_1,\CCZ{S}X_2),$$
\cE{where $(X_1,X_2)$ is a bivariate Gaussian random vector with correlation coefficient $\rho$ and $N(0,1)$ components
being independent of  $\CCZ{S^2}$ which has the} generalized inverse Gaussian distribution 
i.e.,
\BQNY
\pk{\CCZ{S}>x}&\sim &\frac{\left(\frac{\alpha^2}{\delta^2}\right)^{\frac{\lambda}{2}}}
{2\cdot \mathbf{K}_{\lambda}(\sqrt{\delta^2\alpha^2})}\frac{2}{\alpha^2}
x^{2\lambda-2}\exp(-\alpha^2x^2/2)\\
&=&c(\lambda,\delta^2,\alpha^2)\frac{2}{\alpha^2}x^{2\lambda-2}\exp(-\alpha^2x^2/2),
\EQNY
where $\mathbf{K}_{\lambda}$ denotes the modified Bessel function of the third kind (see Abramowitz and Stegun (1965), p. 355-494), \CCZ{$\alpha>0,\delta>0$ and $\lambda \inr$}.\\
Schlueter and Fischer (2012) derived the weak tail dependence coefficient of EGHD by complex calculations. Now using Theorem \ref{ap1}, we immediately obtain the tail asymptotic behaviour and weak tail dependence coefficient for EGHD risks. Indeed, \cE{if  $\CCZ{(Y_1,Y_2)}$ is an EGHD bivariate} random vector defined as above, then Theorem \ref{ap1} yields
\BQN\label{thm}
\pk{\CCZ{Y_1>x,Y_2>x}}\sim\frac{c(\lambda,\delta^2,\alpha^2)}{\sqrt{2\pi}}\frac{(1+\rho)^{3/2}}{(1-\rho)^{1/2}}
\alpha^{-\lambda-\frac{3}{2}}c_{\rho}^{\frac{2\lambda+1}{4}}
x^{\frac{2\lambda-3}{2}}\EXP{-\alpha\sqrt{c_{\rho}}x}
\EQN
and the weak tail dependence coefficient is
$$\chi=2\left(\frac{1+\rho}{2}\right)^{\frac{1}{2}}-1.$$
Note that \eqref{thm} is claimed (but the formula there is not correct) in \cE{Theorem 3} of the aforementioned paper.

\COM{
Our third application is discussed the tail asymptotic of the aggregated risk
$$\pk{Y_1+Y_2>y}, \qquad y \to \IF,$$
where $Y_i, i=1,2$ have tail asymptotics behaviour as
$$\pk{Y_i> y} \sim g_i(\ln y)\EXP{-L (\ln y)^p},$$
with $g_i$ some regularly varying function at infinity and ultimately monotone, further $Y_i$ satisfy the conditions \eqref{depB} and \eqref{depB2}. In view of Theorem \ref{th1},
we have
\BQNY
\pk{Y_1+Y_2>y}&=&\pk{e^{Y_1}e^{Y_2}>e^y}\\
&\sim & D \Bigl(\frac{2\pi p_{2}L_{2}}{p_{1}+p_{2}}\Bigr)
^{1/2}A^{p_{2}/2+q _{2}-q _{1}}g_1(\hat{z}_y^{-1}e^y)g_2(\hat{z}_y)\EXP{\frac{2p_{2} q
_{1}+2p_{1}q _{2}+p_{1}p_{2}}{2(p_{1}+p_{2})}y}
\exp\left(-Be^{\frac{p_1p_2}{p_1+p_2}y}\right),
\EQNY
as $x\to \IF$, where $\hat{z}_y=A\EXP{\frac{p_1}{p_1+p_2}y}$ and $A$ and $B$ given by \eqref{CA}.
}

\section{Proofs}

\prooftheo{th1} \CZC{First 
by} \eqref{depB} we have (recall $w_x=x^{\frac{p_1}{p_1+p_2}}$)
\cW{
\BQNY
\pk{Z> x} =\overline{H}(x)
&=& \cE{\int^{\infty}_{0} \ZW{c(x,y)}\overline{F}_1\left(\frac{x}{y}\right)\,  \d F_2(y)}\\
&\geq&\int_{w_x}^{\infty}c(x,y)\overline{F}_1\left(\frac{x}{y}\right)\,\d F_2(y)\\
&\geq& K_1\int_{w_x}^{\infty}x^{\beta_1}\overline{F}_1\left(\frac{x}{y}\right) \, \d F_2(y)\\
&\geq&\ K_1x^{\beta_1}\CCZ{\overline{F}_1\left(xw_x^{-1}\right)\overline{F}_2(w_x)}.
\EQNY
\COM{
or by \eqref{depB1}
\BQNY
\pk{Z> x} =\overline{H}(x)
&\geq& \pk{X_{1}>x^{\frac{p_{2}}{p_{1}+p_{2}}},  X_{2}>x^{\frac{p_{1}}{p_{1}+p_{2}}}}\\
&\geq& c(x^{\frac{p_{2}}{p_{1}+p_{2}}},x^{\frac{p_{1}}{p_{1}+p_{2}}})\pk{ X_{1}>x^{\frac{p_{2}}{p_{1}+p_{2}}}}
\pk{ X_{2}>x^{\frac{p_{1}}{p_{1}+p_{2}}}}  \\
&\geq& (1-\epsilon)^2
Kx^{\frac{q_1p_{2}+q_2p_1}{p_{1}+p_{2}}}
g_1(x^{\frac{p_2}{p_1+p_2}})g_2(x^{\frac{p_1}{p_1+p_2}})
\exp \left( -(L_{1}+L_{2})x^{\frac{p_{1}p_{2}}{
p_{1}+p_{2}}}\right).
\EQNY}}
%
By the assumptions for some small $a_{1}>0$ we obtain 
\BQNY
\int_0^{a_1 w_x }\overline{F}_1\left(\frac{x}{y}\right)c(x,y)\, \d F_2(y)&\leq& K_2  x^{\beta_2}\int_0^{a_1 w_x} \overline{F}_1\left(\frac{x}{y}\right)\, \d F_2(y)\\
&\leq&K_2  x^{\beta_2} \overline{F}_{1}\left( a_1^{-1}xw_x^{-1}\right)
=o(\overline{H}(x)).
\EQNY
Similarly, for some large $a_{2}>0$ we obtain
\BQNY
\int_{a_2 w_x }^\IF\overline{F}_1\left(\frac{x}{y}\right)c(x,y)\, \d F_2(y)
&\leq& K_2 x^{\beta_2} \overline{F}_{2}\left( a_{2}w_x\right)
=o(\overline{H}(x)).
\EQNY
{Next, by Lemma A.5 in Tang and Tsitsiashvili (2004) 
we can assume that without loss of generality that $F_2$ is absolutely continuous and
therefore we take simply $F_2(x)=1-g_2(x)\exp(-L_2x^{p_2}), x>0$.} Further, in view of Theorem 1.3.3 of Bingham et al.\ (1987) we can assume that $g_2$ is a normalised slowly varying function with derivative $g_2'$. Consequently, for all large $x$ \BQNY
\pk{Z> x} &\sim& \int_{a_1 w_x}^{a_2 w_x} \overline F_1(x/y)c(x,y) \, \d F_2(y)\\
&\sim & -D x^{q_1} \int_{a_1 w_x}^{a_2 w_x} y^{q_2-q_1} \overline F_1(x/y) \, \d \Bigl(g_2(y)\exp(-L_2y^{p_2})\Bigr)\\
&=& D x^{q_1}L_2p_2 \int_{a_1 w_x}^{a_2 w_x} y^{q_2-q_1+p_2-1}
\overline F_1(x/y)\CZC{g_2(y)}\exp(-L_2y^{p_2})\left[1-\CCC{\frac{g_2^{\prime}(y)}{g_2(y)L_2p_2}}y^{1-p_2}\right] \, \d y \\
&\sim& D L_2p_2 x^{q_1}\int_{a_1 w_x}^{a_2 w_x} y^{q_2-q_1+p_2-1}
g_1(x/y)g_2(y)\exp(-L_1(x/y)^{p_1}-L_2y^{p_2}) \, \d y.
\EQNY
We write further
$$I_1(x)+I_2(x)+I_3(x)=
\cE{\Biggl(} \int_{a_{1}w_x}^{(1+\varepsilon)^{-\frac{1}{p_2}}z_x}
+\int_{(1+\varepsilon)^{-\frac{1}{p_2}}z_x}^{(1+\varepsilon)^{\frac{1}{p_2}}z_x}
+\int_{(1+\varepsilon)^{\frac{1}{p_2}}z_x}^{a_{2}w_x}\cE{\Biggr)}\CCZ{ y^{q_2-q_1+p_2-1}
g_1(x/y)g_2(y)\exp(-L_1(x/y)^{p_1}-L_2y^{p_2}) \, \d y},
 $$
where \CCZ{$\varepsilon>0$}, $z_x=Aw_x$ and $ A $ \pE{is} given by \eqref{CA}.
Note that the function $\psi(y)=L_1\left(x/y\right)^{p_1}+L_2y^{p_2}$ decreases when $0<y\leq z_x$ and increases when $y\geq z_x$.
\cE{As \czw{in} Liu and Tang (2010), we obtain}
\begin{eqnarray*}
I_{1}(x)
\leq 
\EXP{ -\left(
\left( 1+\varepsilon \right) ^{\frac{p_{1}}{p_{2}}}L_{1}A^{-p_{1}}+(1+
\varepsilon )^{-1}L_{2}A^{p_{2}}\right) x^{\frac{p_{1}p_{2}}{p_{1}+p_{2}}
}}
\int_{a_{1}w_x}^{(1+\varepsilon)^{-\frac{1}{p_2}}z_x}
g_1\left(\frac{x}{y}\right)g_2(y)y^{q_2-q_1+p_2-1}\, \d y
\end{eqnarray*}
and
\begin{eqnarray*}
I_{3}(x)\leq 
\EXP{ - \left(
(1+\varepsilon )^{-\frac{p_{1}}{p_{2}}}L_{1}A^{-p_{1}}+\left( 1+\varepsilon
\right) L_{2}A^{p_{2}}\right) x^{\frac{p_{1}p_{2}}{p_{1}+p_{2}}}}
\int_{(1+\varepsilon)^{\frac{1}{p_2}}z_x}^{a_{2}w_x}
g_1\left(\frac{x}{y}\right)g_2(y)y^{q_2-q_1+p_2-1}\, \d y.
\end{eqnarray*}
\pE{Next, we have}
\begin{eqnarray*}
&&I_{2}(x) \\
&\geq & 
\left( \int_{(1+\varepsilon /2)^{-\frac{1}{p_{2}}}z_x}^{z_x}+
\int_{z_x}^{(1+\varepsilon /2)^{\frac{1}{p_{2}}}z_x}\right) y^{q_2-q_1+p_2-1}g_1\left(\frac{x}{y}\right)g_2(y)
\EXP{ -\left(L_{1}\fracl{x}{y} ^{p_{1}}+L_{2}y^{p_{2}}\right)}\, \d y \\
&\geq& 
\EXP{ - \left(
\left( 1+\varepsilon /2\right) ^{\frac{p_{1}}{p_{2}}}L_{1}A^{-p_{1}}+\left(
1+\varepsilon /2\right) ^{-1}L_{2}A^{p_{2}}\right)x^{\frac{p_{1}p_{2}}{%
p_{1}+p_{2}}}}
\int_{(1+\varepsilon /2)^{-\frac{1}{p_{2}}}z_x}^{z_x}
y^{q_2-q_1+p_2-1}g_1\left(\frac{x}{y}\right)g_2(y)\, \d y\\
&&+ 
\EXP{ - \left( \left(
1+\varepsilon /2\right) ^{-\frac{p_{1}}{p_{2}}}L_{1}A^{-p_{1}}+\left(
1+\varepsilon /2\right) L_{2}A^{p_{2}}\right) x^{\frac{p_{1}p_{2}}{%
p_{1}+p_{2}}}} \int_{z_x}^{(1+\varepsilon /2)^{\frac{1}{p_{2}}}z_x}
y^{q_2-q_1+p_2-1}g_1\left(\frac{x}{y}\right)g_2(y) \, \d y.
\end{eqnarray*}%
\pE{Consequently, as in Liu and Tang (2010), using Taylor's expansion we obtain}  $I_{1}(x)=o(I_{2}(x))$ and $I_{3}(x)=o(I_{2}(x))$, implying thus
for all $\varepsilon >0$
\begin{eqnarray*}
\overline{H}(x)\sim \CCZ{DL_2p_2x^{q_1}}I_{2}(x)= DL_2p_2x^{q_1}\int_{(1+\varepsilon )^{-\frac{1}{
p_{2}}}z_x}^{(1+\varepsilon )^{\frac{1}{p_{2}}
}z_x} g_1\fracl{x}{y}g_2(y)y^{q_2-q_1+p_2-1}
\EXP{ -L_1 \fracl{x}{y} ^{p_{1}}-L_2 y^{p_{2}}}\, \d y.
\end{eqnarray*}
Since \cE{$g_1(\cdot),g_2(\cdot)$ are ultimately \cQ{monotone}, assume without loss of generality that
they are both ultimately increasing. Hence} for $y \in \left[(1+\varepsilon )^{-\frac{1}{
p_{2}}}z_x,(1+\varepsilon )^{\frac{1}{p_{2}}
}z_x \right]$ we have
\begin{eqnarray}
g_1\left((1+\varepsilon )^{-\frac{1}{p_{2}}
}z_x^{-1}x\right)g_2\left((1+\varepsilon )^{-\frac{1}{p_{2}}
}z_x\right)
\leq g_1\fracl{x}{y}g_2(y)
\leq g_1\left((1+\varepsilon )^{\frac{1}{p_{2}}
}z_x^{-1}x\right)g_2\left((1+\varepsilon )^{\frac{1}{p_{2}}
}z_x\right). \label{uf}
\end{eqnarray}
\cE{By letting} $\varepsilon \to 0$
and using the Laplace approximation we obtain
\begin{eqnarray*}
\overline{H}(x)\sim DL_2p_2x^{q_1}I_2(x)
\sim D\sqrt{\frac{2\pi p_2L_2}{p_1+p_2}}A^{\frac{p_2}{2}+q_2-q_1}x^{\frac{2p_1q_2+2p_2q_1+p_1p_2}{2(p_1+p_2)}}
g_1(z_x^{-1}x)g_2(z_x)
\EXP{-(L_1A^{-p_1}+L_2A^{p_2})w_x^{p_2}},
\end{eqnarray*}
and thus the proof is complete.  \QED

\COM{
and $a$ sufficiently small
\BQNY
\pk{X_1X_2> x} &\sim& \int_{a w_x}^{a^{-1} w_x} \overline F_1(x/y)c(x,y) \, d F_2(y)\\
&\sim & -C_2D x^{q_1} \int_{a w_x}^{a^{-1} w_x} y^{q_2-q_1} \overline F_1(x/y) \, d \Bigl(y^{\alpha_2}\exp(-L_2y^{p_2})\Bigr)\\
&=& C_2D x^{q_1}L_2p_2 \int_{a w_x}^{a^{-1} w_x} y^{q_2-q_1+\alpha_2+p_2-1}
\overline F_1(x/y)\exp(-L_2y^{p_2})\left[1-\frac{\alpha_2}{L_2p_2}y^{-p_2}\right] \, d y\\
&\sim& C_1C_2D L_2p_2 x^{q_1+\alpha_1}\int_{a w_x}^{a^{-1} w_x} y^{q_2-q_1+\alpha_2+p_2-\alpha_1-1}
\exp(-L_1(x/y)^{p_1}-L_2y^{p_2}) \, d y\\
\EQNY
Define the function
$\varphi(y)=L_1(x/y)^{p_1}+L_2y^{p_2}$. Let $\varphi^{\prime}(y)=0$ we get extremum
$z_x=Aw_x$ with $A=[(p_1L_1)/(p_2L_2)]^{1/(p_1+p_2)}$.
Since $\varphi^{\prime\prime}(z_x)>0$, $z_x$ is minimum of $\varphi(y)$ in $[a w_x,a^{-1} w_x]$. Set
$\psi(y)=L_{1}A^{-p_{1}}y^{-p_{1}}+L_{2}A^{p_{2}}y^{p_{2}}$ and
$\widehat \alpha_i=\alpha_i+ q_i,i=1,2$,
applying the Laplace approximation 
we obtain
\BQNY
&&\int_{a w_x}^{a^{-1} w_x} y^{\widehat \alpha_2+p_2-\widehat \alpha_1-1}\exp(-L_1(x/y)^{p_1}-L_2y^{p_2}) \, d y\\
&=&z_x^{\widehat \alpha_2+p_2-\widehat \alpha_1}
\int_{a A^{-1}}^{a^{-1} A^{-1}}y^{\widehat \alpha_2+p_2-\widehat \alpha_1-1}\exp \left(-x^{\frac{p_{1}p_{2}}{p_{1}+p_{2}}%
}[L_{1}A^{-p_{1}}y^{-p_{1}}+L_{2}A^{p_{2}}y^{p_{2}}]\right)  \, d y \\
&=&(1+o(1))z_x^{\widehat \alpha_2+p_2-\widehat \alpha_1}\int_{a A^{-1}-1}^{a^{-1} A^{-1}-1}
\exp\left(-x^{\frac{p_{1}p_{2}}{p_{1}+p_{2}}}
\left(\psi(1)+\frac{\psi^{\prime\prime}(1)}{2}y^2\right)\right) \, d y\\
&=&(1+o(1))z_x^{\widehat \alpha_2+p_2-\widehat \alpha_1}\exp\left(-x^{\frac{p_{1}p_{2}}{p_{1}+p_{2}}}\psi(1)\right)
\int_{a A^{-1}-1}^{a^{-1} A^{-1}-1}
\exp\left(-\frac{1}{2}x^{\frac{p_{1}p_{2}}{p_{1}+p_{2}}}\psi^{\prime\prime}(1)y^2\right) \, d y.\\
\EQNY
Consequently, we obtain the following asymptotic expansion
\begin{eqnarray*}
\pk{X_1X_2> x}\sim \sqrt{2\pi} C_{1}C_{2}D L_{2} p_{2}x^{\widehat \alpha
_{1}}z_{x}^{\widehat \alpha_2+p_2-\widehat \alpha_1}
\exp\left(-x^{\frac{p_{1}p_{2}}{p_{1}+p_{2}}}\psi(1)\right)
\left(\frac{x^{-\frac{p_1p_2}{2(p_1+p_2)}}}{\sqrt{\psi^{\prime\prime}(1)}}
\right),
\end{eqnarray*}
thus the result follows.
}

\proofkorr{co1}
\COM{
Arbitrarily choose $\varepsilon >0$. Since
\begin{eqnarray*}
\overline{H}(x)=\pk{X_1X_2>x} &\geq& \pk{ X_{1}>x^{\frac{p_{2}}{p_{1}+p_{2}}}}
\pk{ X_{2}>x^{\frac{p_{1}}{p_{1}+p_{2}}}}  \\
&\geq& (1-\varepsilon) g_1(x^{\frac{p_2}{p_1p_2}})g_2(x^{\frac{p_1}{p_1+p_2}})
\EXP{ -(L_{1}+L_{2})x^{\frac{p_{1}p_{2}}{
p_{1}+p_{2}}}}
\end{eqnarray*}
it holds for some small $v_{1}>0$ that
\begin{eqnarray*}
\int_0^{a_1x^{\frac{p_1}{p_1+p_2}}}\pk{X_1>\frac{x}{y}}\,dF_2(y)&\leq&
\overline{F_{1}}\left( v_{1}^{-1}x^{\frac{p_{2}}{p_{1}+p_{2}}}\right)\\
&\leq& (1+\varepsilon)g_1(a_1^{-1}x^{\frac{p_2}{p_1+p_2}})
\exp(-L_1a_1^{-p_1}x^{\frac{p_1p_2}{p_1+p_2}})\\
&=&o(\overline{H}(x)).
\end{eqnarray*}
Likewise, it holds for some large $v_{2}>0$ that
\begin{eqnarray*}
\int_{a_2x^{\frac{p_1}{p_1+p_2}}}^{\infty}\pk{X_1>\frac{x}{y}}\,dF_2(y)&\leq&
\overline{F_{2}}\left( v_{2}x^{\frac{p_{1}}{p_{1}+p_{2}}}\right) \\
&\leq& (1+\varepsilon)g_2(a_2x^{\frac{p_1}{p_1+p_2}})
\exp(-L_2a_2^{p_2}x^{\frac{p_1p_2}{p_1+p_2}})\\
&=&o(
\overline{H}(x)).
\end{eqnarray*}
Therefore, we have
\begin{eqnarray*}
\overline{H}(x)=\int_{0}^{\infty }\overline{F_{1}}\fracl{x}{y} \,d
F_{2}(y)&\sim& \int_{v_{1}x^{\frac{p_{1}}{p_{1}+p_{2}}}}^{v_{2}x^{\frac{p_{1}
}{p_{1}+p_{2}}}}\overline{F_{1}}\fracl{x}{y} \,d F_{2}(y)\\
&\sim& L_2p_2\int_{v_{1}x^{\frac{p_{1}}{p_{1}+p_{2}}}}^{v_{2}x^{\frac{p_{1}}{p_{1}+p_{2}}}}
g_1\left(\frac{x}{y}\right)g_2(y)y^{p_2-1}\EXP{-L_1\left(\frac{x}{y}\right)^{p_1}-L_2y^{p_2}}\, dy.
\end{eqnarray*}
Split the integral on the right-hand side above into three parts as
$$I_1(x)+I_2(x)+I_3(x)=
\int_{v_{1}x^{\frac{p_{1}}{p_{1}+p_{2}}}}^{(1+\varepsilon)^{-\frac{1}{p_2}}z_x}
+\int_{(1+\varepsilon)^{-\frac{1}{p_2}}z_x}^{(1+\varepsilon)^{\frac{1}{p_2}}z_x}
+\int_{(1+\varepsilon)^{\frac{1}{p_2}}z_x}^{v_{2}x^{\frac{p_{1}}{p_{1}+p_{2}}}}$$
where $z_x=Ax^{\frac{p_1}{p_1+p_2}}$ and $ A =[(p_1L_1)/(p_2 L_2)]^{1/(p_1+p_2)}$.
Note that the function $\psi(y)=L_1\left(x/y\right)^{p_1}+L_2y^{p_2}$, decreases when $0<y\leq z_x$ and increases when $y\geq z_x$.
Next, we may write
\begin{eqnarray*}
I_{1}(x)
\leq L_2p_2
\exp \left\{ -\left(
\left( 1+\varepsilon \right) ^{\frac{p_{1}}{p_{2}}}L_{1}A^{-p_{1}}+(1+
\varepsilon )^{-1}L_{2}A^{p_{2}}\right) x^{\frac{p_{1}p_{2}}{p_{1}+p_{2}}
}\right\}
\int_{v_{1}x^{\frac{p_{1}}{p_{1}+p_{2}}}}^{(1+\varepsilon)^{-\frac{1}{p_2}}z_x}
g_1\left(\frac{x}{y}\right)g_2(y)y^{p_2-1}\, dy
\end{eqnarray*}
and
\begin{eqnarray*}
I_{3}(x)\leq L_2p_2
\exp \left\{ - \left(
(1+\varepsilon )^{-\frac{p_{1}}{p_{2}}}L_{1}A^{-p_{1}}+\left( 1+\varepsilon
\right) L_{2}A^{p_{2}}\right) x^{\frac{p_{1}p_{2}}{p_{1}+p_{2}}}\right\}
\int_{(1+\varepsilon)^{\frac{1}{p_2}}z_x}^{v_{2}x^{\frac{p_{1}}{p_{1}+p_{2}}}}
g_1\left(\frac{x}{y}\right)g_2(y)y^{p_2-1}\, dy
\end{eqnarray*}
Furthermore, we have
\begin{eqnarray*}
&&I_{2}(x) \\
&\geq & L_2p_2\left( \int_{(1+\varepsilon /2)^{-\frac{1}{p_{2}}}z_x}^{z_x}+
\int_{z_x}^{(1+\varepsilon /2)^{\frac{1}{p_{2}}}z_x}\right) y^{p_2-1}g_1\left(\frac{x}{y}\right)g_2(y)
\exp \left\{ -\left(L_{1}\fracl{x}{y} ^{p_{1}}+L_{2}y^{p_{2}}\right) \right\}\, dy \\
&\geq& L_2p_2
\exp \left\{ - \left(
\left( 1+\varepsilon /2\right) ^{\frac{p_{1}}{p_{2}}}L_{1}A^{-p_{1}}+\left(
1+\varepsilon /2\right) ^{-1}L_{2}A^{p_{2}}\right)x^{\frac{p_{1}p_{2}}{%
p_{1}+p_{2}}}\right\}
\int_{(1+\varepsilon /2)^{-\frac{1}{p_{2}}}z_x}^{z_x}
y^{p_2-1}g_1\left(\frac{x}{y}\right)g_2(y)\, dy\\
&&+ L_2p_2
\exp \left\{ - \left( \left(
1+\varepsilon /2\right) ^{-\frac{p_{1}}{p_{2}}}L_{1}A^{-p_{1}}+\left(
1+\varepsilon /2\right) L_{2}A^{p_{2}}\right) x^{\frac{p_{1}p_{2}}{%
p_{1}+p_{2}}}\right\} \int_{z_x}^{(1+\varepsilon /2)^{\frac{1}{p_{2}}}z_x}
y^{p_2-1}g_1\left(\frac{x}{y}\right)g_2(y) \, dy.
\end{eqnarray*}%
Using Taylor's expansion we get $I_{1}(x)=o(I_{2}(x))$ and $I_{3}(x)=o(I_{2}(x))$ (see Liu and Tang (2010)).
\COM{
to expand both $\left( 1+\varepsilon /2\right) ^{%
\frac{p_{1}}{p_{2}}}L_{1}A^{-p_{1}}+\left( 1+\varepsilon /2\right)
^{-1}L_{2}A^{p_{2}}$ and $(1+\varepsilon )^{\frac{p_{1}}{p_{2}}%
}L_{1}A^{-p_{1}}+\left( 1+\varepsilon \right) ^{-1}L_{2}A^{p_{2}}$ in $%
\varepsilon $ up to $\varepsilon ^{2}$. Then we find that the coefficients
of the constant terms and the $\varepsilon $ terms are equal, but the
coefficient of the $\varepsilon ^{2}$ term of the former is smaller than the
corresponding coefficient of the latter. Therefore, for all small $%
\varepsilon >0$,
\begin{eqnarray*}
&& \left( \left( 1+\varepsilon /2\right) ^{%
\frac{p_{1}}{p_{2}}}L_{1}A^{-p_{1}}+\left( 1+\varepsilon /2\right)
^{-1}L_{2}A^{p_{2}}\right)  \\
&<& \left( (1+\varepsilon )^{\frac{p_{1}}{%
p_{2}}}L_{1}A^{-p_{1}}+\left( 1+\varepsilon \right)
^{-1}L_{2}A^{p_{2}}\right) ,
\end{eqnarray*}%
which implies that $I_{1}(x)=o(I_{2}(x))$.
Similarly as above for all small
$\varepsilon >0$
\begin{eqnarray*}
&& \left( \left( 1+\varepsilon /2\right) ^{-%
\frac{p_{1}}{p_{2}}}L_{1}A^{-p_{1}}+\left( 1+\varepsilon /2\right)
L_{2}A^{p_{2}}\right)  \\
&<& \left( \left( 1+\varepsilon \right) ^{-%
\frac{p_{1}}{p_{2}}}L_{1}A^{-p_{1}}+(1+\varepsilon )L_{2}A^{p_{2}}\right) ,
\end{eqnarray*}%
which implies that $I_{3}(x)=o(I_{2}(x))$.}
Hence, for all $\varepsilon >0$ we have
\BQNY
\overline{H}(x)\sim I_{2}(x)&=& L_2p_2\int_{(1+\varepsilon )^{-\frac{1}{
p_{2}}}z_x}^{(1+\varepsilon )^{\frac{1}{p_{2}}
}z_x} g_1\fracl{x}{y}g_2(y)y^{p_2-1}
\exp \left\{ -L_1 \fracl{x}{y} ^{p_{1}}-L_2 y^{p_{2}}\right\}\, dy. \label{FI}
\EQNY
Since $g_i(\cdot)$ are monotone, without loss of generality, assume that
$g_i(\cdot)$ are increase, $i=1,2$, then for $y \in [(1+\varepsilon )^{-\frac{1}{
p_{2}}}z_x,(1+\varepsilon )^{\frac{1}{p_{2}}
}z_x ]$ we have
\begin{eqnarray*}
g_1((1+\varepsilon )^{-\frac{1}{p_{2}}
}z_x^{-1}x)g_2((1+\varepsilon )^{-\frac{1}{p_{2}}
}z_x)
\leq g_1\fracl{x}{y}g_2(y)
\leq g_1((1+\varepsilon )^{\frac{1}{p_{2}}
}z_x^{-1}x)g_2((1+\varepsilon )^{\frac{1}{p_{2}}
}z_x).
\end{eqnarray*}
Let $\varepsilon \to 0,$
using the Laplace approximation we get
\begin{eqnarray*}
\overline{H}(x)\sim I_2(x) &\sim&
L_2p_2g_1(z_x^{-1}x)g_2(z_x)
\int_{(1+\varepsilon )^{-\frac{1}{
p_{2}}}z_x}^{(1+\varepsilon )^{\frac{1}{p_{2}}
}z_x} y^{p_2-1}
\exp \left\{ -L_1 \fracl{x}{y} ^{p_{1}}-L_2 y^{p_{2}}\right\}\, dy\\
&\sim&L_2p_2g_1(z_x^{-1}x)g_2(z_x)
\left(z_x\right)^{p_2}
\int_{(1+\varepsilon )^{-\frac{1}{
p_{2}}}}^{(1+\varepsilon )^{\frac{1}{p_{2}}
}} y^{p_2-1}
\exp \left\{ -x^{\frac{p_1p_2}{p_1+p_2}}(L_1A^{-p_1}y^{-p_1}+L_2A^{p_2}y^{p_2})\right\}\, dy\\
&\sim& \sqrt{\frac{2\pi p_2L_2}{p_1+p_2}}A^{\frac{p_2}{2}}x^{\frac{p_1p_2}{2(p_1+p_2)}}
g_1(z_x^{-1}x)g_2(z_x)
\EXP{-(L_1A^{-p_1}+L_2A^{p_2})x^{\frac{p_1p_2}{p_1+p_2}}}.
\end{eqnarray*}
}
\COM{
Notice that
\BQNY
h(x)&=&\int^\infty_0 h_1\left(\frac{x}{y}\right)\frac{1}{y}\mathrm{d}F_2(y)\\
&=&(\int_0^{y_0}+\int_{y_0}^{\infty})h_1\left(\frac{x}{y}\right)\frac{1}{y}\mathrm{d}F_2(y)\\
&=&J_1+J_2,
\EQNY
and
\BQNY
\limsup_{x \to \infty}\frac{J_1}{J_2}&=&\limsup_{x \to \infty}
\frac{\int^{y_0}_{0}h_1\left(\frac{x}{y}\right)\frac{1}{y}\mathrm{d}F_2(y)}
{\int_{y_0}^{\infty} h_1\left(\frac{x}{y}\right)\frac{1}{y}\mathrm{d}F_2(y)}\\
&\leq&\frac{1}{\int_{y^*}^{\infty}\frac{1}{y}\mathrm{d}F_2(y)}
\int^{y_0}_{0}\limsup_{x \to \infty}\frac{h_1(\frac{x}{y})}{h_1(\frac{x}{y^*})}\frac{1}{y}\mathrm{d}F_2(y)\\
&=&0.
\EQNY
Hence
\BQNY
h(x)\sim \int_{y_0}^{\infty} h_1\left(\frac{x}{y}\right)\frac{1}{y}\mathrm{d}F_2(y).
\EQNY
Let $X^*$ be a positive random variable with distribution function $F^*(x)$
satisfies
\begin{eqnarray*}
 \overline{F^*}(x) \sim  x^{p_1} g_1(x) \exp(- L_1 x^{p_1}),
\end{eqnarray*}
then arbitrarily choose $\epsilon >0$, we have
\begin{eqnarray*}
h(x)&=&\int^\infty_0h_1\left(\frac{x}{y}\right)\frac{1}{y}\mathrm{d}F_2(y)\\
&\geq&(1-\epsilon^2)\int^\infty_0  L_1p_1 \left(\frac{x}{y}\right)^{p_1-1}
\overline{F}_1\left(\frac{x}{y}\right)\frac{1}{y}\mathrm{d}F_2(y)\\
&\geq&(1-\epsilon^2)(1-\epsilon)L_1p_1x^{-1}\int^\infty_0  \left(\frac{x}{y}\right)^{p_1}g_1\left(\frac{x}{y}\right)
\exp\left(-L_1\left(\frac{x}{y}\right)^{p_1}\right)\mathrm{d}F_2(y)\\
&\geq&(1-\epsilon)^2L_1p_1x^{-1}\int^\infty_0  \overline{F^*}\left(\frac{x}{y}\right)\mathrm{d}F_2(y)\\
&\geq&(1-\epsilon)^2L_1p_1x^{-1}\pk{X^*>x^{\frac{p_2}{p_1+p_2}}}\pk{X_2>x^{\frac{p_1}{p_1+p_2}}}\\
&\geq&(1-\epsilon)^4L_1p_1x^{p_1-1}g_1(x^{\frac{p_2}{p_1+p_2}})g_2(x^{\frac{p_1}{p_1+p_2}})
\exp \left\{ -(L_{1}+L_{2})x^{\frac{p_{1}p_{2}}{p_{1}+p_{2}}}\right\},
\end{eqnarray*}}
The tail asymptotic of the distribution of $Z$ \cE{follows easily, therefore we show next} the tail asymptotic of the pdf $h$ of $Z$.
For all $x$ large \cE{and $\epsilon\in (0,1)$}, since  $h_1$ is ultimately decreasing
\BQNY
h(x)&=&\int^\infty_0h_1\left(\frac{x}{y}\right)\frac{1}{y}\, \d F_2(y)\\
&\geq&\int_{w_x}^{2w_x}h_1\left(\frac{x}{y}\right)\frac{1}{y}\, \d F_2(y)\\
&\geq&2^{-1}w_x^{-1}h_1(xw_x^{-1})[\overline{F}_2(w_x)-\overline{F}_2(2w_x)]\\
&\geq&\cE{(1-\epsilon)}2^{-1}L_1 p_1 x^{p_1-\frac{p_1^{\czw{2}}}{p_1+p_2}-1}
g_1(xw_x^{-1})g_2(w_x)
\EXP{ -(L_{1}+L_{2})w_x^{p_2}}\\
&&\times\left(1-\frac{g_2(2w_x)}{g_2(w_x)}\EXP{-L_2(2^{p_2}-1)w_x^{p_2}}\right).
\EQNY
Let $X^*$ be a positive random variable with distribution function $F^*$ which satisfies
\begin{eqnarray*}
 \overline{F}^*(x) \sim  x^{p_1} g_1(x) \exp(- L_1 x^{p_1}).
\end{eqnarray*}
\cE{For} \czw{some} $a_{1}>0$ small enough we have
\begin{eqnarray*}
\int_0^{a_1w_x}h_1\left(\frac{x}{y}\right)\frac{1}{y}\,  \d F_2(y)
&\leq&(1+\epsilon)L_1p_1x^{-1}\int_0^{a_1w_x}
\overline{F}^*\left(\frac{x}{y}\right)\, \d F_2(y)\\
&\leq&(1+\epsilon)L_1p_1x^{-1}\overline{F}^*\left( a_{1}^{-1}w_x^{-1}x\right)
=o(h(x)).
\end{eqnarray*}
\cE{Similarly,} for some large $a_{2}>0$, \CCZ{since $h_1$ is bounded, there \pE{exists a} positive constant $M$ such that}
\BQNY
\int_{a_2w_x}^{\infty}h_1\left(\frac{x}{y}\right)\frac{1}{y} \, \d F_2(y)
\leq  Ma_2^{-1}w_x^{-1}g_2(a_{2}w_x)
\EXP{-L_2a_{2}^{p_2}w_x^{p_2}}=o(h(x)).
\EQNY
\cE{Consequently,} with the same arguments as in the proof of Theorem \ref{th1}  we obtain
\BQNY
h(x)
&\sim& \int_{a_{1}w_x}^{a_{2}w_x}h_1\left( \frac{x}{y}\right)\frac{1}{y}  \, \d F_{2}(y)\\
&\sim& L_1L_2p_1p_2x^{p_1-1}\int_{a_{1}w_x}^{a_{2}w_x}
y^{p_2-p_1-1}g_1\left(\frac{x}{y}\right)g_2(y)\EXP{-L_1\left(\frac{x}{y}\right)^{p_1}-L_2y^{p_2}}\, \d y,
\EQNY
hence the proof follows easily. \QED

\COM{
\prooftheo{th3}
Arbitrarily choose $\epsilon >0$. Since
\begin{eqnarray*}
\overline{H}(x) &\geq& \pk{ X_{1}>x^{\frac{p_{2}}{p_{1}+p_{2}}}}
\pk{ X_{2}>x^{\frac{p_{1}}{p_{1}+p_{2}}}}  \\
&\geq& (1-\epsilon)^2C_1C_2x^{\frac{p_1\alpha_2+p_2\alpha_1}{p_1+p_2}}
\exp \left\{ -(L_{1}+L_{2})x^{\frac{p_{1}p_{2}}{
p_{1}+p_{2}}}\right\}
\end{eqnarray*}
it holds for some large $v>0$ that
\begin{eqnarray*}
\int_0^{v^{-1}x^{\frac{p_1}{p_1+p_2}}}\overline{F}_1\left(\frac{x}{y}\right) dF_2(y)
\leq
\overline{F_{1}}\left( vx^{\frac{p_{2}}{p_{1}+p_{2}}}\right)
\leq (1+\epsilon)C_1v^{\alpha_1} x^{\frac{\alpha_1p_{2}}{p_{1}+p_{2}}}
\exp(-L_1v^{p_1}x^{\frac{p_1p_2}{p_1+p_2}})
=o\left(\EXP{-(L_1+L_2)vx^{\frac{p_1p_2}{p_1+p_2}}}\right).
\end{eqnarray*}
and
\begin{eqnarray*}
\int_{vx^{\frac{p_1}{p_1+p_2}}}^{\infty}\overline{F}_1\left(\frac{x}{y}\right) dF_2(y)
\leq
\overline{F_{2}}\left( vx^{\frac{p_{1}}{p_{1}+p_{2}}}\right)
\leq (1+\epsilon)C_2v^{\alpha_2} x^{\frac{\alpha_2p_{1}}{p_{1}+p_{2}}}
\exp(-L_2v^{p_2}x^{\frac{p_1p_2}{p_1+p_2}})
=o\left(\EXP{-(L_1+L_2)vx^{\frac{p_1p_2}{p_1+p_2}}}\right).
\end{eqnarray*}
Therefore,
\begin{equation*}
\overline{H}(x)=\int_{0}^{\infty }\overline{F_{1}}\left( x/y\right) \mathrm{d%
}F_{2}(y)= \int_{v^{-1}x^{\frac{p_{1}}{p_{1}+p_{2}}}}^{vx^{\frac{p_{1}%
}{p_{1}+p_{2}}}}\overline{F_{1}}\left( x/y\right) \mathrm{d}F_{2}(y)
+o\left(\EXP{-(L_1+L_2)vx^{\frac{p_1p_2}{p_1+p_2}}}\right).
\end{equation*}
Since $X_2$ possess a positive pdf $h_2(x) =(1+o(1))\pk{X_2> x}  L_2p_2 x^{p_2-1}$, we have
\begin{eqnarray*}
\overline{H}(x)&=&(1+o(1))C_1C_2L_2p_2x^{\alpha_1}
\int_{v^{-1}x^{\frac{p_{1}}{p_{1}+p_{2}}}}^{vx^{\frac{p_{1}}{p_{1}+p_{2}}}}
y^{p_2+\alpha_2-\alpha_1-1}\exp\left(-L_1\left(x/y\right)^{p_1}-L_2y^{p_2}\right)\mathrm{d}y
+o\left(\EXP{-(L_1+L_2)vx^{\frac{p_1p_2}{p_1+p_2}}}\right).
\end{eqnarray*}
Since the function
$\varphi(y)=L_1(x/y)^{p_1}+L_2y^{p_2}$ obtains its minimum in $[a_1x^{\frac{p_1}{
p_1+p_2}},a_2x^{\frac{p_1}{p_1+p_{2}}
}]$
at $Ax^{\frac{p_1}{p_1+p_2}}$, set $z_x=Ax^{\frac{p_1}{p_1+p_2}}$ and
$\psi
(y)=L_{1}A^{-p_{1}}y^{-p_{1}}+L_{2}A^{p_{2}}y^{p_{2}}$,
applying the Laplace approximation
we obtain
\begin{eqnarray*}
&&\int_{v_{1}x^{\frac{p_{1}}{p_{1}+p_{2}}}}^{v_{2}x^{\frac{p_{1}}{p_{1}+p_{2}}}}
y^{p_2+\alpha_2-\alpha_1-1}\exp\left(-L_1\left(x/y\right)^{p_1}-L_2y^{p_2}\right)\mathrm{d}y\\
&=&z_x^{p_2+\alpha_2-\alpha_1}
\int_{a_1A^{-1}}^{a_2A^{-1}}y^{p_2+\alpha_2-\alpha_1-1}\exp \left(-x^{\frac{p_{1}p_{2}}{p_{1}+p_{2}}%
}[L_{1}A^{-p_{1}}y^{-p_{1}}+L_{2}A^{p_{2}}y^{p_{2}}]\right) \mathrm{d}y \\
&=&(1+o(1))z_{x}^{p_{2}+\alpha _{2}-\alpha _{1}}
\exp\left(-x^{\frac{p_{1}p_{2}}{p_{1}+p_{2}}}\psi(1)\right)
\left(\frac{x^{-\frac{p_1p_2}{2(p_1+p_2)}}}{\sqrt{\psi^{\prime\prime}(1)}}+
\sum^{\infty}_{n=1}\lambda_n x^{-\frac{(2n+1)p_{1}p_{2}}{2(p_{1}+p_{2})}}
\right)
\end{eqnarray*}
where the coefficients $\lambda_n$ are polynomials of the higher derivatives $\psi^{(k)}(1),k\geq2$.
Hence
\begin{eqnarray*}
\overline{H}(x)&=&(1+o(1))C_1C_2L_2p_2x^{\alpha_1}
z_{x}^{p_{2}+\alpha _{2}-\alpha _{1}}
\exp\left(-x^{\frac{p_{1}p_{2}}{p_{1}+p_{2}}}\psi(1)\right)
\left(\frac{x^{-\frac{p_1p_2}{2(p_1+p_2)}}}{\sqrt{\psi^{\prime\prime}(1)}}+
\sum^{\infty}_{n=1}\lambda_n x^{-\frac{(2n+1)p_{1}p_{2}}{2(p_{1}+p_{2})}}\right)\\
&&+o\left(\EXP{-(L_1+L_2)vx^{\frac{p_1p_2}{p_1+p_2}}}\right)\\
&=&(1+o(1))C_1C_2L_2p_2x^{\alpha_1}
z_{x}^{p_{2}+\alpha _{2}-\alpha _{1}}
\exp\left(-x^{\frac{p_{1}p_{2}}{p_{1}+p_{2}}}\psi(1)\right)
\left(\frac{x^{-\frac{p_1p_2}{2(p_1+p_2)}}}{\sqrt{\psi^{\prime\prime}(1)}}+
\sum^{\infty}_{n=1}\lambda_n x^{-\frac{(2n+1)p_{1}p_{2}}{2(p_{1}+p_{2})}}\right),
\end{eqnarray*}
hence he claim follows. \QED
}

\prooftheo{ap1}  In view of  \eqref{elliptical}
\BQNY
\pk{SX_1>x,SX_2>x}&\sim&
\frac{1}{2\pi}\frac{(1-\rho^2)^{\frac{3}{2}}}{(1-\rho)^2}\frac{p_1+p_2}{p_1p_2}(L_1A^{-p_1}+L_2A^{p_2})^{-1}\\
&& \times \left(\frac{2}{1+\rho}\right)^{-\frac{p_1p_2}{2(p_1+p_2)}+1}
x^{-\frac{p_1p_2}{p_1+p_2}}\pk{SR>\sqrt{\frac{2}{1+\rho}}x}
\EQNY
and thus by \eqref{eq FGM} the first \cE{claim follows}.
Since $Z$ is in the Gumbel max-domain of attraction with scaling function $w(\cdot)$ given by
$$ w(x)= (L_1A^{-p_1}+L_2A^{p_2})\frac{p_1p_2}{p_1+p_2}\CCZ{x}^{\frac{p_1p_2}{p_1+p_2}-1},$$
then  by Theorem 2.1 of Hashorva (2010)  for any $c>0$ we have
$$
\lim_{x \to \infty}\frac{ w(cx)}{ w(x)}
=c^{\frac{p_1p_2}{p_1+p_2}-1},
$$
hence the weak tail dependence coefficient is
\BQNY \chi=2\cdot\left(\frac{1+\rho}{2}\right)^{\frac{p_1p_2}{2(p_1+p_2)}}-1
\EQNY
establishing thus the proof. \QED

\COM{
\section{Appendix} We show next the claim in \eqref{Gauss}.
For any $x$ positive with $\Phi$ the distribution function of a $N(0,1)$ random variable we have
\BQNY
\pk{Z>x}&=&\overline{H}(x)
=\int_{\pE{-\IF}}^{\infty} \left(1-\Phi\left(\frac{x/y-\rho y}{\sqrt{1-\rho^2}}\right)\right)\, \d \Phi(y).
\EQNY
\cE{Hence, for some positive constant $C$ \pE{(which may change below from line to line)} and all large $x$
\BQNY
\cE{\overline{H}(x)\geq \pE{\Bigl(}1- \Phi( \sqrt{x}) \Bigr) \Bigl( 1- \Phi(  \sqrt{ x(1- \rho)/(1+\rho)} )\Bigr) \geq \frac{C}{x}\EXP{-x/(1+ \rho)}.}
\EQNY
}
Next, for some small $a_1>0$ 
\BQNY
\int_{0}^{a_1\sqrt{x}}\left(1-\Phi\left(\frac{x/y-\rho y}{\sqrt{1-\rho^2}}\right)\right)\,\d\Phi(y)
&\leq&1-\Phi\left(\sqrt{x}\frac{a_1^{-1}-\rho a_1}{\sqrt{1-\rho^2}}\right)\\
&\leq& \pE{C x^{-1/2}}\EXP{-\frac{x}{2(1-\rho^2)}(a_1^{-1}-a_1\rho)^2}\\
&=&o(\overline{H}(x)).
\EQNY
Similarly, for some large $a_2>0$ we obtain
\BQNY
\int_{a_2\sqrt{x}}^{\infty}\left(1-\Phi\left(\frac{x/y-\rho y}{\sqrt{1-\rho^2}}\right)\right)\,\d\Phi(y)
\leq1-\Phi(a_2\sqrt{x})
\leq \pE{C x^{-1/2}}\EXP{-\frac{a_2^2}{2}x} =o(\overline{H}(x)).
\EQNY
\pE{Consequently, we may further write}
\BQNY
\overline{H}(x)
&\sim& 2\int_{a_1\sqrt{x}}^{a_2\sqrt{x}}\left(1-\Phi\left(\frac{x/y-\rho y}{\sqrt{1-\rho^2}}\right)\right)\,\d\Phi(y)\\
&\sim&   \sqrt{\frac{2(1-\rho^2)}{\pi}}\int_{a_1\sqrt{x}}^{a_2\sqrt{x}}\left(\frac{x}{y}-\rho y\right)^{-1}
\EXP{-\frac{1}{2(1-\rho^2)}\left(\frac{x}{y}-\rho y\right)^2}\, \d\Phi(y)\\
&\sim&\frac{\sqrt{1-\rho^2}}{\CCC{\pi}}\EXP{\frac{\rho}{1-\rho^2}x}
\int_{a_1\sqrt{x}}^{a_2\sqrt{x}}\left(\frac{x}{y}-\rho y\right)^{-1}
\EXP{-\frac{1}{2(1-\rho^2)}\left(\frac{x^2}{y^2}+y^2\right)} \, \d y\\
&=&\frac{\sqrt{1-\rho^2}}{ \CCC{\pi} }\EXP{\frac{\rho}{1-\rho^2}x}
\int_{a_1}^{a_2}y\left(1-\rho y^2\right)^{-1}
\EXP{-\frac{1}{2(1-\rho^2)}x\left(y^{-2}+y^2\right)} \, \d y\\
&\sim&\CCC{\frac{1}{\pi}\sqrt{\frac{1+\rho}{1-\rho}}\EXP{-\frac{x}{1+\rho}}
\int_{a_1-1}^{a_2-1}\EXP{-\frac{2x}{1-\rho^2}y^2}\,\d y}\\
&\sim&{\frac{1+\rho}{\sqrt{2\pi x}}}\EXP{-\frac{x}{1+\rho}},
\EQNY
and thus the proof of \eqref{Gauss} is complete.
}

\textbf{Acknowledgments}: We would like to thank both referees for constructive suggestions that improved our paper significantly.
Z. Weng kindly acknowledges the full financial support by the Swiss National Science Foundation Grant 200021-134785.

\end{document}